\documentclass[letterpaper, 10 pt, conference]{ieeeconf} 

\IEEEoverridecommandlockouts 
\usepackage{graphicx} 
\usepackage{amsmath}

\usepackage{amsthm}
\usepackage{amssymb}
\usepackage{geometry}
\usepackage{cite}
\usepackage{graphicx}
\usepackage{subcaption}
\let\labelindent\relax
\usepackage{enumitem}
\usepackage[hidelinks]{hyperref}
\usepackage{cellspace}
\usepackage{xcolor}
\usepackage[linesnumbered,lined,boxed,longend,commentsnumbered,ruled]{algorithm2e}
\usepackage{algpseudocode}

\def\BibTeX{{\rm B\kern-.05em{\sc i\kern-.025em b}\kern-.08em
    T\kern-.1667em\lower.7ex\hbox{E}\kern-.125emX}}
\newcommand{\norm}[1]{\left \Vert #1 \right \Vert}
\newcommand{\euclidnorm}[1]{\left \vert #1 \right \vert}
\newcommand{\R}{\ensuremath{\mathbb{R}}}
\newcommand{\Z}{\ensuremath{\mathbb{Z}}}

\theoremstyle{definition}

\title{\LARGE \bf Policy iteration for discrete-time systems with discounted costs: stability and near-optimality guarantees}
\author{Jonathan de Brusse$^{1}$, Mathieu Granzotto$^{2}$, Romain Postoyan$^{1}$ and  Dragan Ne\v si\'c$^{2}$%
\thanks{* This work was funded by Lorraine Université d’Excellence LUE and supported by the ANR under grant OLYMPIA ANR-23-CE48-0006.}
\thanks{$^{1}$ J. de Brusse and R. Postoyan are with the Université de Lorraine, CNRS, CRAN, F-54000, Nancy, France. (emails: jonathan.de-brusse@univ-lorraine.fr, romain.postoyan@univ-lorraine.fr).}%
\thanks{$^{2}$ D. Ne\v si\'c and M. Granzotto are with the Department of Electrical and Electronic Engineering, University of Melbourne, Parkville, VIC 3010, Australia (emails: dnesic@unimelb.edu.au, mathieu.granzotto@unimelb.edu.au).}%
}

\date{October 2023}

\begin{document}

\maketitle

\begin{abstract}
Given a discounted cost, we study deterministic discrete-time systems whose inputs are generated by policy iteration (PI). We provide novel near-optimality and stability properties, while allowing for non stabilizing initial policies. That is, we first give novel bounds on the mismatch between the value function generated by PI and the optimal value function, which are less conservative in general than those encountered in the dynamic programming literature for the considered class of systems.  Then, we show that the system in closed-loop with policies generated by PI are stabilizing under mild conditions, after a finite (and known) number of iterations.
\end{abstract}

\section{Introduction}
Dynamic programming provides powerful methods to generate near-optimal inputs for general dynamical systems and cost functions \cite{Bertsekas-book12(adp)}. To make the best out of dynamic programming algorithms in a control engineering context, it is often essential to endow the obtained closed-loop system with stability guarantees. Various results exist in the literature ensuring stability properties for systems controlled by dynamic programming, both in continuous-time and discrete-time, see, e.g., \cite{Granzotto-et-al-journal-submission, Postoyan-et-al-tac(optimal), Zhang-et-al-book2012, Jiang-Jiang-book2017, Vrabie, heydari-acc2016 }. The vast majority of these works focus on \textit{undiscounted} cost functions. However, discounted costs are ubiquitous in dynamic programming and reinforcement learning \cite{Bertsekas-book12(adp), SuttonBarto:98}, because of the favorable properties the discount endows to Bellman operators, like contractivity, see, e.g., \cite{Bertsekas-book12(adp)}. We may also consider discounted costs because of the problem at hand or when no admissible policy (a policy  with a finite bound) for the undiscounted is known, while an admissible policy  for the discounted costs is known. It is therefore important to provide stability guarantees for systems controlled by dynamic programming algorithms with discounted costs.

In \cite{Granzotto-et-al-journal-submission, Granzotto-et-al-learnconf2021} stability results are provided for discrete-time systems controlled by value iteration with discounted costs. Results for discounted policy iteration (PI) are only available for linear systems with quadratic costs (LQ)\cite{Lamperski_2020} where the discount factor is not fixed but increases with the number of iterations, as far as we know. In this work, we consider general deterministic discrete-time systems and costs with fixed discount factors. Our main goal is to establish stability properties when the inputs are generated by PI as well as novel near-optimality bounds. 

We make several assumptions for this purpose. We first assume that an optimal sequence of inputs exists for any initial state and is stabilizing, which is very natural in the context of this work. We also assume that PI is recursively feasible in the sense that the optimization problem solved at each iteration is guaranteed to always admit a solution as customary in the literature \cite{ liu-wei-tnnls13}; if this is not the case we can resort to the modification of PI as advocated in \cite{Granzotto-arxiv22} and our results apply. On the other hand, the initial policy for PI is required to give a bounded finite cost, which, does not mean that it is necessarily stabilizing because of the discount factor. This allows to relax the conditions of the related  literature, see, e.g. \cite{Vrabie, heydari-acc2016,Granzotto-arxiv22}, which require the initial policy to be stabilizing. This is one additional possible reason to consider discounted costs, i.e., to remove the need for an initial stabilizing policy. Finally, the system needs to satisfy a detectability property with respect to the stage cost, which is also very natural as we aim to establish stability properties.

We use a generic measuring function to define stability as in e.g., \cite{Grimm-et-al-tac2005, Postoyan-et-al-tac(optimal)}, which is useful to study the stability of the origin and of more general attractors in a unified way. By exploiting the assumed stability property verified by the system in closed-loop with optimal controllers, novel near-optimality properties are deduced. Indeed, less conservative bounds on the mismatch between the value function generated by PI and the optimal value function compared to \cite{Bertsekas-book12(adp),Munos-03} are obtained by exploiting the properties of the class of systems. Afterwards, we present a Lyapunov analysis, which allows to establish that, after a sufficient number of iterations, PI generates stabilizing policies. This is so even in the absence of stability properties of  the initial policy as already mentioned. In particular, we show that the closed-loop system controlled by PI enjoys a semiglobal practical stability property where the adjustable parameter is the number of iterations. We provide an explicit relationship between the number of iterations required, the set of initial conditions and the guaranteed ultimate bound. By strengthening the assumptions, a global exponential stability property is derived and easy-to-compute lower bound  on the discount factor and the number of iteration are given. We illustrate our results by means of two examples: the LQ problem and a nonholonomic integrator.

The rest of the paper is organized as follows. Preliminaries are recalled in Section II. The problem is formally stated in Section III. The standing assumptions are presented in Section IV. The main results are given in Section V. Examples are presented in Section VI before concluding in Section VII. Long proofs are postponed to the appendix for the sake of readability.

\section{Preliminaries}

\subsection{Notation}
Let $\R$ be the set of real numbers, $\R_{\ge0}:=[0,+\infty)$, $\Z_{\ge0}:=\{0,1,2,...\}$ and $\Z_{>0}:=\{1,2,...\}$. We consider $\mathcal{K}$, $\mathcal{K}_{\infty}$ and $\mathcal{KL}$ functions as defined in \cite[Section 3.5]{Goebel-Sanfelice-Teel-book}. The identity map from $\R_{\ge0}$ to $\R_{\ge0}$ is denoted by $\mathbb{I}$. Let $f : \R_{\ge0} \rightarrow \R_{\ge 0}$, we use $f^{(k)}$ for the composition of function $f$ to itself $k$ times, where $k \in \Z_{\ge 0}$ and $f^{(0)} := \mathbb{I}$. We use $\left \lceil{\cdot} \right\rceil$ to denote the ceil function. The Euclidean norm of a vector $x\in \R^{n}$ with $n\in \Z_{>0}$ is denoted by $\euclidnorm{x}$ and the distance of $x$ to a non-empty set $\mathcal{A}\subseteq\R^{n}$ by $\euclidnorm{x}_{\mathcal{A}}:=\inf\{\euclidnorm{x-y}: y\in \mathcal{A}\}$. For any $M\in \R^{n \times m}$ with $n,m\in \Z_{>0}$, $\norm{M}$ is the spectral norm of the matrix $M$, i.e., $\norm{M}=\sqrt{\rho(M^\top M)}$, where $\rho(M^\top M)$ is the spectral radius of the matrix $M^\top M$. If $M$ is a symmetric matrix let $\lambda_{\min}(M)$ and $\lambda_{\max}(M)$ denote its minimal and maximal eigenvalues, respectively. Given a set-valued map $S: \R^{n}\rightrightarrows\R^{m}$, a selection of $S$ is a single-valued mapping $s:\text{dom}\, S \rightarrow \R^{m}$ such that $s(x)\in S(x)$ for any $x\in \text{dom} \,S$, we write $s\in S$ to denote a selection $s$ of $S$ for the sake of convenience. Finally, for an infinite sequence $\boldsymbol{u}=(u(0),u(1),...)$ where $u(0),u(1),...\in \R^{m}$ with $m\in \Z_{>0}$, $\boldsymbol{u}|_k$ stands for the truncation of $\boldsymbol{u}$ to its first $k\in \Z_{\ge 0}$ steps, i.e., $\boldsymbol{u}|_k=\big(u(0),...,u(k-1)\big)$ and we use the convention $\boldsymbol{u}|_0=\emptyset$.

\subsection{Plant Model and Cost Function}
We consider nonlinear deterministic discrete-time systems given by  
\begin{equation}
    x(k+1)=f(x(k),u(k)), \qquad \forall k \in \Z_{\ge0},
\label{plant}
\end{equation}
where $x \in \R^{n_x}$ is the state, $u\in \mathcal{U}(x) \subseteq \R^{n_u}$ is the control input, $\mathcal{U}(x)$ is the non-empty set of \textit{admissible} inputs at state $x\in \R^{n_x}$, and $n_x,n_u\in \Z_{>0}$. Ideally, we wish to find, for any given $x\in \R^{n_x}$, an infinite-length sequence of admissible inputs $\boldsymbol{u}=\big(u(0),u(1),...\big)$ that minimizes the discounted infinite-horizon cost
\begin{equation}
    J_{\gamma}(x,\boldsymbol{u}):=\displaystyle \sum_{k=0}^{\infty} \gamma^k \ell\big(\phi(k,x,\boldsymbol{u}|_k),u(k)\big),
    \label{total cost}
\end{equation}
where $\gamma\in (0,1)$ is a discount factor, $\ell: \R^{n_x} \times \R^{n_u} \rightarrow \R_{\ge0}$ is a non-negative stage cost and $\phi(k,x,\boldsymbol{u}|_k)$ is the solution to \eqref{plant} at time $k\in \Z_{\ge0}$, initialized at $x$, with inputs given by $\boldsymbol{u}|_k$ and we use the convention $\phi(0,x,\boldsymbol{u}|_0)=x$. We assume that for any $x\in \R^{n_x}$, there exists a sequence of admissible inputs minimizing $J_{\gamma}(x,\cdot)$, i.e., 
\begin{equation}
    V^{\star}_{\gamma}(x):=\underset{\boldsymbol{u}}{\min} \, J_{\gamma}(x,\boldsymbol{u})<+\infty, \qquad \forall x\in \R^{n_x},
    \label{Minimal cost}
\end{equation}
as formalized in Section III-A. As a consequence, Bellman equation becomes
\begin{equation}
    \label{Bellman equation}
    V^{\star}_{\gamma}(x)= \underset{u \in \mathcal{U}(x)}{\min}\big\{\ell(x,u)+\gamma V^{\star}_{\gamma}\big(f(x,u)\big)\big\} \qquad \forall x\in \R^{n_x}.
\end{equation}
We can therefore define the non-empty set of optimal inputs for any state $x \in \R^{n_x}$ as
\begin{equation}
    H^{\star}_{\gamma}(x):=\underset{u\in \mathcal{U}(x)}{\text{argmin}} \{ \ell(x,u)+\gamma V^{\star}_{\gamma}\big(f(x,u)\big)\}.
    \label{Set of optimal policies}
\end{equation}
Given \eqref{Set of optimal policies}, the closed loop system \eqref{plant} with optimal controller is given by
\begin{equation}
\label{Optimal plant}
    x(k+1)\in f\big(x(k),H^{\star}_{\gamma}(x(k))\big)=: F^{\star}_{\gamma}\big(x(k)\big) \quad \forall k\in \Z_{\ge 0}.
\end{equation}
As \eqref{Set of optimal policies} is a set-valued map, there may be non-unique optimal inputs at some state and, as a consequence, system \eqref{Optimal plant} is a difference inclusion in general. For the sake of convenience, solutions to system \eqref{Optimal plant} at time $k\in \Z_{\ge0}$ are denoted as $\phi^{\star}_{\gamma}(k,x)$ when initialized at $x\in \R^{n_x}$.

Computing $H^{\star}_{\gamma}$ in \eqref{Set of optimal policies} for the general dynamics in $\eqref{plant}$ and cost function \eqref{total cost} is notoriously hard. Dynamic programming provides algorithms to iteratively obtain feedback laws, which instead generate policies that converge to an optimal one. We focus on policy iteration (PI) in this work, which we recall in the next section. Before that, we introduce some notation. Given a policy $h:\R^{n_x}\rightarrow \R^{n_u}$ that is admissible, i.e., $h\in \mathcal{U}$, we denote the solution to system \eqref{plant} in closed-loop with feedback law $h$ at time $k\in \Z_{\ge 0}$ with initial condition $x$ as $\phi(k,x,h)$. Likewise $J_{\gamma}(x,h)$ is the cost induced by $h$ at initial state $x$,  i.e., $J_{\gamma}(x,h):=\displaystyle \sum_{k=0}^{\infty} \gamma^k \ell\big(\phi(k,x,h),h\big(\phi(k,x,h)\big)\big)$.

\section{Problem statement}
We recall PI in this section and we state the objectives of this work.
\subsection{Policy Iteration}
PI is given in Algorithm 1. Given $\gamma \in (0,1)$ and an initial admissible policy $h^{0}$, PI generates at each iteration $i\in \Z_{\ge0}$ a policy $h^{i+1}_{\gamma}$ via the so-called improvement step in \eqref{PI.2}. Policy $h^{i+1}_{\gamma}$ is an arbitrary selection of $H^{i+1}_{\gamma}$ in \eqref{PI.2} where $H^{i+1}_{\gamma}$ may be set-valued. We then evaluate the cost induced by $h^{i+1}_{\gamma}$, namely $V^{i+1}_{\gamma}(x)=J_{\gamma}(x,h^{i+1}_{\gamma})$ for any $x\in \R^{n_x}$, at the evaluation step in \eqref{PI.3}. By doing so repeatedly, $V^{i}_{\gamma}$ converges to the optimal value function $V^{\infty}_{\gamma}=V^{\star}_{\gamma}$ under mild conditions, see,  e.g., \cite{Bertsekas-book12(adp)}.

 It is implicitly assumed here that the optimization problem defined in \eqref{PI.2} always admits a solution, i.e., $H^i_{\gamma}(x)$ is non-empty for any $x\in \R^{n_x}$ at any iteration $i\in \Z_{>0}$. We say in this case that Algorithm 1 is \textit{recursively feasible}. We will go back to this point in Section IV-B.
 
\begin{algorithm}[]
\SetAlgoNoLine
\KwIn{$f$ in \eqref{plant}, $\ell$ in \eqref{total cost}, $\gamma \in (0,1)$, initial policy $h^0\in \mathcal{U}$}
\KwOut{Policy $h^{\infty}_{\gamma}$, cost $V^{\infty}_{\gamma}$}
    \textbf{Initial evaluation step:} for all $x\in \R^{n_x}$
    \begin{equation}
        \label{PI.1}
        \tag{PI.1}
        V^{0}_{\gamma}(x):=J_{\gamma}(x,h^0).
    \end{equation}\\
    \For{$i \in \mathbb{Z}_{\ge 0}$}{
    \text{\textbf{Policy improvement step:} for all $x\in \R^{n_x}$}\begin{equation}
        \label{PI.2}
        \tag{PI.2}
        H^{i+1}_{\gamma}(x):=\underset{u \in \mathcal{U}(x)}{\text{argmin}} \{\ell(x,u)+\gamma V^{i}_{\gamma}(f(x,u))\}.
    \end{equation}\\
    \textbf{Select } $h^{i+1}_{\gamma} \in H^{i+1}_{\gamma}$.\\
    \textbf{Policy evaluation step:} for all $x\in \R^{n_x}$,
    \begin{equation}
        \label{PI.3}
        \tag{PI.3}
        V^{i+1}_{\gamma}(x):=J_{\gamma}(x,h^{i+1}_{\gamma}).
    \end{equation}}
 \textbf{return } $h^{\infty}_{\gamma} \in H^{\infty}$ and $V^{\infty}_{\gamma}$.
\caption{Policy Iteration}
\end{algorithm}

\subsection{Objectives}

Our objective is to give conditions under which PI generates stabilizing policies after a finite number of iterations. We also aim at providing near-optimality guarantees for PI. In particular, we will see that the bound on $V^i_{\gamma} - V^{\star}_{\gamma}$ we provide significantly differ and improve those encountered in the dynamic programming literature \cite{Bertsekas-book12(adp), Munos-03} for the considered class of systems.

We need to make several assumptions to meet these objectives, which are presented in the next section.

\section{Standing Assumptions}

\subsection{Existence of an Optimal Sequence}
As mentioned in Section II-B, we assume that for any $x\in \R^{n_x}$, there exists (at least) one infinite-length sequence of admissible inputs minimizing \eqref{total cost}.

\textit{Standing Assumption 1 (SA1):} For any $x\in \R^{n_x}$ and any $\gamma \in (0,1)$,  there exists an optimal sequence of admissible inputs $\boldsymbol{u}^{\star}_{\gamma}(x)$ such that $V^{\star}_{\gamma}(x)=J_{\gamma}\big(x,\boldsymbol{u}_{\gamma}^{\star}(x)\big)<+\infty$ and for any infinite-length sequence of admissible inputs $\boldsymbol{u}$, $V^{\star}_{\gamma}(x)\leq J_{\gamma}(x,\boldsymbol{u})$. \hfill $\Box$

Condition on system \eqref{plant} and cost function \eqref{total cost} ensuring SA1 are available in \cite{Keerthi-Gilbert-tac85} for instance. SA1 ensures the existence of the optimal value function $V^{\star}_{\gamma}$ given in \eqref{Minimal cost} as well as the non-emptiness of $H^{\star}_{\gamma}(x)$ in \eqref{Set of optimal policies} for any $x \in \R^{n_x}$ and $\gamma\in(0,1)$. SA1 is very reasonable in the context of this work as we aim to use PI to generate policies that converge to an optimal one, which therefore needs to exist.

\subsection{Recursive Feasibility of PI}

We proceed as is often done in the literature, see, e.g., \cite{heydari-acc2016, Vrabie, liu-wei-tnnls13}, and assume Algorithm 1 is recursively feasible, in the sense that the set $H^{i}_{\gamma}(x)$ is non-empty at any iteration $i\in \Z_{>0}$ and for any $x\in \R^{n_x}$, which is equivalent to say that the optimization problem in \eqref{PI.2} admits a solution for any $x\in \R^{n_x}$ at any iteration $i\in \Z_{>0}$. 

\textit{Standing Assumption 2 (SA2):} For any $i\in \Z_{>0}$ and $x\in \R^{n_x}$, the set-valued map $H^{i}_{\gamma}(x)$ is non-empty. \hfill$\Box$

SA2 ensures the recursive feasibility of Algorithm 1 by allowing the selection at each iteration of a new policy. As explained in the introduction, if SA2 does not hold, we can use the modified version of PI presented in \cite[Section IV]{Granzotto-arxiv22} and the forthcoming results apply mutatis mutandis.

\subsection{Detectability}
To define stability, we use a continuous and radially unbounded function\footnote{In the sense that for any $\Delta>0$, $\{x\in \R^{n_x}: \sigma(x)\leq \Delta\}$ is compact.} $\sigma: \R^{n_x} \rightarrow \R_{\ge0}$ that serves as a state measure relating the distance of the state to a given attractor where $\sigma$ vanishes. As explained in \cite{Grimm-et-al-tac2005, Postoyan-et-al-tac(optimal)}, $\sigma$ can be defined as $\left \vert \cdot \right \vert^p$ when studying the stability of the origin, or as $\left \vert \cdot \right \vert^{p}_{\mathcal{A}}$, with $p\in \Z_{>0}$ when studying the stability of non-empty set $\mathcal{A}\subseteq \R^{n_x}$.  We make the next detectability assumption on system \eqref{plant} and stage cost $\ell$, which is inspired by \cite{Grimm-et-al-tac2005}.

\textit{Standing Assumption 3 (SA3):} There exist a continuous function $W:\R^{n_x}\rightarrow\R_{\ge0}$, $\alpha_{W}\in \mathcal{K}_{\infty}$ and $\overline{\alpha}_{W}: \R_{\ge0}\rightarrow \R_{\ge0}$ continuous, nondecreasing and zero at zero, such that, for any $x\in \R^{n_x}$ and $u\in \mathcal{U}(x)$,
\begin{equation}
\label{Detectability assumption}
\begin{split}
    W(x)&\leq \overline{\alpha}_{W}\big(\sigma(x)\big)\\
    W\big(f(x,u)\big)-W(x)&\leq -\alpha_{W}\big(\sigma(x)\big)+\ell(x,u).
\end{split}
\end{equation} \hfill $\Box$

SA3 is a detectability property of system \eqref{plant} with respect to $\sigma$ when considering $\ell$ as an output, which is very natural as this captures the fact, that by minimizing $\ell$ along the solutions to \eqref{plant}, desirable stability properties should follow. SA3 is consistent with the literature on LQ \cite{Anderson-Moore-book}, the link between \eqref{Detectability assumption} and detectability of linear time-invariant systems being established in  \cite[Lemma 4]{Postoyan-et-al-tac(optimal)}

\textit{Remark 1:} A more general detectability property is considered in \cite{Grimm-et-al-tac2005, Postoyan-et-al-tac(optimal)} where the second line of \eqref{Detectability assumption} is given by $W(f(x,u))-W(x)\leq -\alpha_{W}(\sigma(x))+\chi(\ell(x,u))$ for any $x\in \R^{n_x}$, $u\in \R^{n_u}$ and with $\chi\in \mathcal{K}_{\infty}$. We plan to investigate this generalization in future work. \hfill $\Box$

\subsection{Initial Policy}

We also make the next assumption on the cost given by the initial policy. 

\textit{Standing Assumption 4 (SA4):} Let $h^0\in \mathcal{U}$ be known and such that there exist $\gamma_0\in(0,1]$ and $\overline{\alpha}_{V}: \R_{\ge0}\times(0,\gamma_0) \rightarrow \R_{\ge 0}$ of class $\mathcal{K}_{\infty}$ in its first argument verifying for any $x\in \R^{n_x}$ and any $\gamma\in (0,\gamma_0)$,
\begin{equation}
    V^{0}_{\gamma}(x)=J_{\gamma}(x,h^{0})\leq \overline{\alpha}_{V}\big(\sigma(x),\gamma\big).
\end{equation} \hfill $\Box$

SA4 requires to know an initial policy $h^0$ such that there exists a range of values for $\gamma$, namely $(0,\gamma_0)$ with $\gamma_0$ sufficiently small, such that the initial cost $V^{0}_{\gamma}(x)$ is finite for every $x \in \R^{n_x}$ and is upper-bounded by function $\overline{\alpha}_V$, which is $\mathcal{K}_{\infty}$ in $\sigma(x)$ and depends on $\gamma$. It is important to note that SA4 may hold even if $h^{0}$ is not stabilizing as illustrated in Section VI and exemplified below. The next lemma gives a sufficient condition to ensure SA4.

\textit{Lemma 1:} Consider system \eqref{plant} and suppose there exist $M, a>0$, $\chi\in \mathcal{K}_{\infty}$ and an admissible policy $h \in \mathcal U$ such that for any $x\in \R^{n_x}$ and any $k\in \Z_{\ge0}$, $\ell\big(\phi(k,x,h),h(\phi(k,x,h))\big)\leq Ma^k\chi(\sigma(x))$. Then SA4 is verified with $h^0 =h$, $\gamma_0=\min\big\{1,\tfrac{1}{a}\big\}$ and $\overline{\alpha}_V(s,\gamma)=\tfrac{M}{1-a\gamma}\chi(s)$ for any $s\in \R_{\ge0}$ and any $\gamma \in (0,\gamma_0)$. \hfill $\Box$\\
\begin{proof} Let $h^0=h$ and $\gamma\in (0,\gamma_0)$  with $h$ and $\gamma_0$ as in Lemma 1. For any $x\in \R^{n_x}$, 
\begin{align}
    V^{0}_{\gamma}(x) &=\displaystyle \sum_{k=0}^{\infty} \gamma^k \ell\big(\phi(k,x,h^0),h(\phi(k,x,h^0))\big)\nonumber\\
    &\leq \displaystyle M\chi(\sigma(x))\sum_{k=0}^{\infty} (a\gamma)^k.
\end{align}
As $\gamma<\gamma_0\leq \tfrac{1}{a}$, we finally have $V^{0}_{\gamma}(x) \leq \tfrac{M}{1-a\gamma}\chi(\sigma(x))$.
 Thus SA4 holds with $\overline{\alpha}_V(\cdot,\gamma):=\tfrac{M}{1-a\gamma}\chi$ and this concludes the proof.
\end{proof}

Lemma 1 shows that, if the stage cost along the solution to \eqref{plant} with policy $h^0$ is upper-bounded at any time-step $k\in \Z_{\ge 0}$ by a term $Ma^k\sigma(x)$, we can determine explicitly $\gamma_0$ and $\overline{\alpha}_V(\cdot,\gamma)$ verifying SA4. Note that $a$ in Lemma 1 may be strictly bigger than 1, which implies that $h^0$ may not be stabilizing. 

\subsection{Stability with Optimal Sequence}
As we want to eventually obtain stabilizing policies using PI, we need to assume that optimal policies are stabilizing. The next assumption together with SA1 and SA3 indeed guarantee that the closed loop system \eqref{plant} with optimal controller, i.e., system \eqref{Optimal plant}, verifies a $\mathcal{KL}$-stability property with respect to $\sigma$ as established in Proposition 1 below. 

\textit{Standing Assumption 5 (SA5):} The following holds.
\begin{enumerate}[label=(\roman*)]
    \item There exists $\overline{\alpha}_{V^{\star}} \in \mathcal{K}_{\infty}$ such that for any $\gamma \in (0,\gamma_0)$ with $\gamma_0$ in SA4, for any $x \in \R^{n_x}$, $V^{\star}_{\gamma}(x)\leq \overline{\alpha}_{V^{\star}}\big(\sigma(x)\big)$ where $V^{\star}_{\gamma}$ is the optimal value function in \eqref{Minimal cost}.
    \item  There exists $\gamma^{\star}\in(0,\gamma_0)$ such that
\begin{equation}
\label{Assumption 3}
    (1-\gamma^{\star})\overline{\alpha}_{V^{\star}}(s)\leq\alpha_W(s), \qquad \forall s\in \R_{>0},
\end{equation}
with $\alpha_W$ in SA3. \hfill $\Box$
\end{enumerate} 

Item (i) may be obtained by studying the case where $\gamma=\gamma_0$ as $V^{\star}_{\gamma}(x)\leq V^{\star}_{\gamma_0}(x)$ for any $x\in \R^{n_x}$ and any $\gamma \in (0,\gamma_0)$. It is important to note that we do not need to know either $V^{\star}_{\gamma}$ or $V^{\star}_{\gamma_0}$ to check whether item (i) of SA5 holds. Indeed, this condition holds whenever the cost for a known, not necessarily optimal, policy is upper-bounded by some $\overline{\alpha}_{V^{\star}}(\sigma(x))$ for any $x\in \R^{n_x}$. A condition ensuring item (i) is given in \cite[Lemma 1]{Postoyan-et-al-tac(optimal)}. Item (ii) is a technical condition useful to deduce global asymptotic stability properties for system \eqref{Optimal plant}, as formalized next.

\textit{Proposition 1:} For any $\gamma\in(\gamma^{\star},\gamma_0)$, system \eqref{Optimal plant} is \emph{$\mathcal{KL}$-stable with respect to $\sigma$}, i.e., there exists $\beta^{\star}_{\gamma}\in \mathcal{KL}$ such that for any $x\in \R^{n_x}$, any solution $\phi^{\star}_{\gamma}(\cdot,x)$ to system \eqref{Optimal plant} satisfies 
\begin{equation}
    \sigma\big(\phi^{\star}_{\gamma}(k,x)\big)\leq \beta^{\star}_{\gamma}\big(\sigma(x),k\big) \quad \forall k\in \Z_{\ge0}.
\end{equation}
In particular, $\beta^{\star}_{\gamma}:(s,k) \mapsto \underline{\alpha}_{Y^{\star}}^{-1}(\widetilde{\beta}^{\star}_{\gamma}(\overline{\alpha}_{Y^{\star}}(s),k)) \in \mathcal{KL}$ with $\underline{\alpha}_{Y^{\star}}$, $\widetilde{\beta}^{\star}_{\gamma}$ and $\overline{\alpha}_{Y^{\star}}$ in Table 1. \hfill $\Box$

\textit{Remark 2:} Proposition 1 ensures a global asymptotic stability property for system \eqref{Optimal plant}. We will investigate in future work the case where a semiglobal practical stability holds for \eqref{Optimal plant} instead as in \cite{Postoyan-et-al-tac(optimal)}, which will allow us to relax item (ii) of SA5. \hfill $\Box$

Now that all the standing assumptions have been stated, we are ready to present the main results of this work.

\section{Main results}
In this section, we consider system $\eqref{plant}$ whose inputs are generated by PI at iteration $i\in \Z_{\ge 0}$, that is,
\begin{equation}
\label{Set function i}
    x(k+1)\in f\big(x(k), H^{i}_{\gamma}(x(k))\big)=:F^{i}_{\gamma}\big(x(k)\big), \quad \forall i\in \Z_{\ge 0}. 
\end{equation}
 For convenience, solutions to system \eqref{Set function i} are denoted in the sequel as $\phi^{i}_{\gamma}(k,x)$ when initialized at $x\in \R^{n_x}$ for any $k\in \Z_{\ge 0}$.
 
\subsection{Near-optimality}
We first recover the classical result presented in, e.g.,  \cite{Bertsekas-Tsitsiklis-book96}, on the improvement property of the policies generated by PI, whose proof is omitted as it follows similar lines as \cite[Lemma 2]{heydari-acc2016}.

\textit{Lemma 2:} For any $x\in \R^{n_x}$, $i\in \Z_{\ge0}$ and $\gamma\in (0,\gamma_0)$ with $\gamma_0$ in SA4, $V^{i+1}_{\gamma}(x)\leq V^{i}_{\gamma}(x)$. \hfill $\Box$

In the next theorem, we establish a new near-optimality bound for PI with discounted costs.

\textit{Theorem 1:}  For any $i\in \Z_{\ge0}$, $x\in \R^{n_x}$, $\gamma  \in \big(\gamma^{\star},\gamma_0\big)$ and any solution $\phi^{\star}_{\gamma}(\cdot,x)$ to system \eqref{Optimal plant},
\begin{align}
    \label{Convergence PI 1}
    \big(V^{i}_{\gamma}-V^{\star}_{\gamma}\big)(x)& \leq \gamma^i\big(V^{0}_{\gamma}-V^{\star}_{\gamma}\big)(\phi^{\star}_{\gamma}(i,x))\nonumber \\
    &\leq \gamma^i\overline{\alpha}_{V}(\beta^{\star}_{\gamma}(\sigma(x),i),\gamma),
\end{align}
with $\beta^{\star}_{\gamma}\in \mathcal{KL}$ from Proposition 1 and $\overline{\alpha}_{V}$ from SA4.

\hfill $\Box$

Theorem 1 gives us an explicit upper-bound on the term $V^{i}_{\gamma}(x)-V^{\star}_{\gamma}(x)$ for any $x\in \R^{n_x}$ and any iteration $i\in \Z_{\ge0}$. Hence, this bound provides an estimation on how close cost $V^{i}_{\gamma}$ generated at iteration $i\in \Z_{\ge 0}$ by PI is from the ``target'' $V^{\star}_{\gamma}$. As $\beta^{\star}_{\gamma}\in \mathcal{KL}$ for any $\gamma\in(\gamma^{\star},\gamma_0)$, the upper bound in \eqref{Convergence PI 1} converges to 0 as the number of iteration $i$ goes to infinity. Typical near-optimal bounds for PI in the literature are of the form $\tfrac{M_1 \gamma^i}{1-\gamma} + M_2$ with $M_1, M_2 \in \R_{>0}$, see, e.g., \cite{Bertsekas-book12(adp), Munos-03}. The latter bound explodes as $\gamma$ goes to 1 and does not vanish to 0 as $\sigma(x)$ goes to 0, contrary to the bound in Theorem 1. This can be explained by the fact that the bounds of the literature do not exploit the stability properties of the optimal policies. 

We now focus on the stability guarantees that can be deduced thanks to this near-optimality property.

\subsection{Stability}
We first establish a Lyapunov property for system \eqref{Set function i} for any iteration $i\in \Z_{\ge0}$.

\textit{Theorem 2:} There exist $\underline{\alpha}_Y \in \mathcal{K_{\infty}}$,  $\overline{\alpha}_{Y}$, $\alpha_Y: \R_{\ge 0} \times (\gamma^{\star},\gamma_0) \rightarrow \R_{\ge 0}$ of class $\mathcal{K}_{\infty}$ in their first argument such that for any $i\in \Z_{\ge 0}$ there exist $Y^{i}_{\gamma}:\R^{n_x} \rightarrow \R_{\ge0}$ and $\Upsilon^{i}: \R_{\ge0}\times (\gamma^{\star},\gamma_0)\rightarrow \R_{\ge0}$ of class $\mathcal{K}_{\infty}$ in its first argument such that the following holds for any $\gamma\in (\gamma^{\star},\gamma_0)$.
\begin{enumerate}[label=(\roman*)]
    \item For any $x\in \R^{n_x}$, $\underline{\alpha}_{Y}\big(\sigma(x)\big)\leq Y^{i}_{\gamma}(x) \leq \overline{\alpha}_{Y}\big(\sigma(x),\gamma\big)$.
    \item For any $x\in \R^{n_x}$ and $v\in F^{i}_{\gamma}(x)$,\\
    $Y^{i}_{\gamma}(v)-Y^{i}_{\gamma}(x)\leq \tfrac{1}{\gamma} \Big( -\alpha_Y\big(\sigma(x),\gamma\big)+\Upsilon^{i}\big(\sigma(x),\gamma\big)\Big)$,
\end{enumerate}
where $\underline{\alpha}_{Y}$, $\overline{\alpha}_Y,\alpha_Y$, $\Upsilon^i$ and $Y^{i}_{\gamma}$ are defined in Table I.\hfill $\Box$

Item (i) of Theorem 1 means that $Y^{i}_{\gamma}$ is positive definite and radially unbounded with respect to $\sigma$ for any $\gamma\in(\gamma^{\star},\gamma_0)$. Item (ii) of Theorem 2 is a dissipative inequality of system \eqref{Set function i} for which the supply rate consists of a negative term, namely $-\alpha_Y(\cdot,\gamma)$, and a non-negative term $\Upsilon^i(\cdot,\gamma)$ which can be made as small as desired by increasing $i$. The latter property is key for establishing stability properties for system $\eqref{Set function i}$ with $i$ sufficiently large. 

\begin{table}[h]
    \centering
    \begin{tabular}{ll}
    \hline
    $\underline{\alpha}_{Y^{\star}}$, $\underline{\alpha}_{Y}$ & $\alpha_W$\\
    $\alpha_{Y^{\star}}(\cdot,\gamma)$, $\alpha_Y(\cdot,\gamma)$ & $\tfrac{\gamma-\gamma^{\star}}{1-\gamma^{\star}}\alpha_W$\\
    $\overline{\alpha}_{Y^{\star}}$ & $ \overline{\alpha}_{V^{\star}}+\tfrac{1}{\gamma^{\star}}\overline{\alpha}_W$\\
    $\widetilde{\alpha}_{Y^{\star}}(\cdot,\gamma)$ & $\alpha_Y(\cdot,\gamma)\circ \overline{\alpha}_{Y^{\star}}$ \\
    $\widetilde{\beta}^{\star}_{\gamma}(s,k)$ & $\max_{\widehat{s}\in [0,s]}\big(\mathbb{I}-\tfrac{1}{\gamma}\widetilde{\alpha}_{Y^{\star}}(\cdot,\gamma)\big)^{(k)}(\widehat{s})$\\
    $Y^{i}_{\gamma}$ & $V^{i}_{\gamma}+\tfrac{1}{\gamma}W$ \\
    $\overline{\alpha}_Y(\cdot,\gamma)$ & $\overline{\alpha}_{V}(\cdot,\gamma)+\tfrac{1}{\gamma}\overline{\alpha}_W$\\
    $\Upsilon^i$ & $(s_1,s_2)\mapsto (1-s_2)s_2 ^i\overline{\alpha}_{V}\big(\beta^{\star}_{\gamma}(s_1,i)\big)$\\
    \hline
    \end{tabular}
    \caption{Expressions of functions used in Theorems 1 and 2}
    \label{tab:my_label}
\end{table}
 Based on Theorem 2, we establish the next stability property for system \eqref{Optimal plant}.

\textit{Theorem 3}: For any $\gamma \in (\gamma^{\star},\gamma_0)$, there exists $\beta_{\gamma}\in \mathcal{KL}$ such that for any $\delta, \Delta>0$, $i\ge i^{\star}_{\gamma}$ with $i^{\star}_{\gamma}\in \Z_{\ge0}$ verifying 
\begin{equation}
\label{i_star_semi_global}
   i^{\star}_{\gamma}\ge \tfrac{\ln{\big(\tfrac{\alpha_Y(\overline{\alpha}_{Y}^{-1}(\underline{\alpha}_Y(\delta),\gamma),\gamma)}{2(1-\gamma)\overline{\alpha}_V(\beta^{\star}_{\gamma}(\underline{\alpha}_Y^{-1}(\overline{\alpha}_Y(\Delta,\gamma)),0),\gamma)}}\big)}{\ln{(\gamma)}},
\end{equation}
any $x\in \{z\in \R^{n_x}: \sigma(z)\leq \Delta\}$, any solution $\phi^{i}_{\gamma}(\cdot, x)$ to system \eqref{Set function i} satisfies
\begin{equation}
\label{semiglobal stability}
    \sigma\big(\phi^{i}_{\gamma}(k,x)\big)\leq \max\{\beta_{\gamma}\big(\sigma(x),k\big),\delta\}\quad \forall k\in \Z_{\ge0}. 
\end{equation} \hfill $\Box$

Theorem 3 ensures a semiglobal and pratical stability property of \eqref{Set function i} for any $\gamma \in (\gamma^{\star},\gamma_0)$ where the tuning parameter is the number of iterations $i$. In particular, for any set of initial conditions of the form $\{x\in \R^{n_x}: \sigma(x)\leq \Delta\}$ where $\Delta>0$ can be arbitrarily large, for any arbitrarily small $\delta>0$, there exists $i^{\star}_{\gamma} \in \Z_{\ge0}$ verifying \eqref{i_star_semi_global} such that for any iteration $i\ge i^{\star}_{\gamma}$ \eqref{semiglobal stability} holds. By strengthening the conditions of Theorem 3, it is possible to  derive stronger stability guarantees.

\textit{Corollary 1:}  Suppose there exist $\overline{a}_W\ge0$, $a_W,\overline{a}_{V^{\star}}>0$ and $\overline{a}_V: (\gamma^{\star},\gamma_0) \mapsto \R_{>0}$ such that $\overline{\alpha}_V(s,\cdot)=\overline{a}_V(\cdot) s,$ $\overline{\alpha}_{V^{\star}}(s)=\overline{a}_{V^{\star}} s$, $\alpha_W(s)=a_W s$ and $\overline{\alpha}_W(s)\leq \overline{a}_W s$ for any $s\ge0$ with $\gamma^{\star}=\tfrac{\overline{a}_{V^{\star}}-a_W}{\overline{a}_{V^{\star}}}$. Then for any $\gamma \in \big(\gamma^{\star},\gamma_0\big)$ there exists $i^{\star}_{\gamma} \in \mathbb{Z}_{\ge 0}$ verifying 
\begin{equation}
\label{estimate i}
    i^{\star}_{\gamma}\ge\tfrac{\ln{\big(\tfrac{\gamma^{\star}(\gamma-\gamma^{\star})a_W^2}{2\gamma(1-\gamma)^2\overline{a}_V(\gamma)(\gamma^{\star}\overline{a}_{V^{\star}}+\overline{a}_W)}\big)}}{\ln{\big(\gamma - \tfrac{\gamma^{\star}(\gamma-\gamma^{\star})a_W}{(1-\gamma)(\gamma^{\star}\overline{a}_{V^{\star}}+\overline{a}_W)
    }\big)}},
\end{equation}
such that for any $i\ge i^{\star}_{\gamma}$, $x\in \R^{n_x}$ any solution $\phi^{i}_{\gamma}(\cdot,x)$ to system \eqref{Set function i} satisfies 
\begin{equation}
    \sigma\big(\phi^{i}_{\gamma}(k,x)\big)\leq c_1(\gamma)\sigma(x)e^{-c_2(\gamma)k} \quad \forall k\in\Z_{\ge0}
\end{equation} 
with $c_1(\gamma)=\tfrac{\gamma \overline{a}_V(\gamma)+\overline{a}_W}{\gamma a_W}$ and $c_2(\gamma)=-\ln{\big(1-\tfrac{a_W(\gamma-\gamma^{\star})}{2\gamma(1-\gamma)(\overline{a}_V(\gamma)+\tfrac{1}{\gamma}\overline{a}_W)}}\big)>0$. \hfill $\Box$

\indent  Corollary 1 ensures that after a sufficient number of iterations $i^{\star}_{\gamma}$, that we can explicitly estimate using \eqref{estimate i}, a global exponential stability property of \eqref{Set function i} is also verified for any $\gamma \in (\gamma^{\star},\gamma_0)$. 

\textit{Remark 3:} Corollary 1 also ensures a global exponential stability property of the system \eqref{Optimal plant} with a lower bound for $\gamma^{\star}$ less conservative than those presented in \cite[Corollary 2]{Postoyan-et-al-tac(optimal)} and \cite[Lemma 2]{Granzotto-et-al-journal-submission}. Indeed $\gamma^{\star}=\tfrac{\overline{a}_{V^{\star}}-a_W}{\overline{a}_{V^{\star}}}= 1 - \tfrac{a_W}{\overline{a}_{V^{\star}}} \leq 1 - \tfrac{a_W}{\overline{a}_{V^{\star}}+\overline{a}_W}=:\gamma^{\star}_{[6]}$ and as $1>\tfrac{(\overline{a}_{V^{\star}}-a_W)}{\overline{a}_{V^{\star}}}\tfrac{(\overline{a}_{V^{\star}}+a_W)}{\overline{a}_{V^{\star}}}=\tfrac{\overline{a}_{V^{\star}}^2-a_W^2}{\overline{a}_{V^{\star}}^2}$, we have $\gamma^{\star}=\tfrac{\overline{a}_{V^{\star}}-a_W}{\overline{a}_{V^{\star}}}\leq \tfrac{\overline{a}_{V^{\star}}}{\overline{a}_{V^{\star}}+a_W}=:\gamma^{\star}_{[17]}$. \hfill $\Box$

\section{Examples}

We consider two examples, namely the linear quadratic problem and a nonholonomic integrator, for which we show that the standing assumptions hold thereby implying that the results of Section V apply.

\subsection{Linear Quadratic Problem}

We consider the deterministic linear time-invariant system
\begin{equation}
\label{linear plant}
x(k+1)=Ax(k)+Bu(k),
\end{equation}
where $A\in \R^{n_x \times n_x}$, $B\in \R^{n_x \times n_u}$ and $(A,B)$ stabilizable. Let $\sigma(x)=\norm{x}^2$ and $\ell(x,u)=x^{\top}Qx+u^{\top}Ru$ for any $x\in \R^{n_x}$ and $u\in \R^{n_u}$, where $Q=C^{\top}C \in \R^{n_x\times n_x}$ with $(A,C)$ detectable, and $R\in \R^{n_u \times n_u}$ is a symmetric, positive definite matrix. We set $\mathcal{U}(x)=\R^{n_u}$ for any $x\in \R^{n_x}$.

First, as $(A,B)$ is stabilizable and $(A,C)$ is detectable SA1 holds by \cite{Bertsekas-book12(adp)}. In addition, SA3 is verified with $W(x)=x^{\top}S_2x$ for any $x\in \R^{n_x}$, $\alpha_W(s)=\lambda_{\min}(S_1)s$ and $\overline{\alpha}_W(s)=\lambda_{\max}(S_2)s$ for any $s\in \R_{\ge 0}$, where $S_1,S_2$ are symmetric positive definite matrices, which satisfy
\begin{equation}
\label{LMI}
    \begin{pmatrix}
     A^{\top}S_2A-S_2+S_1- Q & A^{\top}S_2B\\
     B^{\top}S_2A & B^{\top}S_2B- R
\end{pmatrix} \leq 0.
\end{equation}
Note that there always exist such matrices $S_1$ and $S_2$ by \cite[Lemma 4]{Postoyan-et-al-tac(optimal)}. Furthermore, the conditions of Lemma 1 are satisfied by taking $h(x)=K_0x$ for any $x\in \R^{n_x}$, with any $K_0 \in \R^{n_u \times n_x}$, $\chi=\mathbb{I}$, $M=\norm{Q+K_0^{\top}RK_0}$ and $a=\norm{A+BK_0}^2$. Hence SA4 is verified with  $\gamma_0 = \min \Big\{1, \tfrac{1}{\norm{A+BK_0}^2}\Big\}$ and $\overline{\alpha}_{V}\big(s,\gamma\big):= \tfrac{\norm{Q+K_0^{\top}RK_0}}{1-\gamma\norm{A+BK_0}^2} s$ for any $s\in \R_{\ge0}$ and $\gamma \in (0,\gamma_0)$. We note that SA4 is verified for any $K_0\in \R^{n_u \times n_x}$, thus even when $A+B K_0$ is not Schur, i.e., even when the initial policy is not stabilizing. Moreover, using the same time-varying change of coordinates as in \cite{Lamperski_2020}, we know by \cite{Vrabie} that SA2 holds, as $\sqrt{\gamma}(A+BK_0)$ is Schur for any $\gamma<\gamma_0$ (which does not mean $A+BK_0$ is Schur obviously). Given $(A,B)$ stabilizable, the optimal value function for any $\gamma\in [0,1]$ and $x\in \R^{n_x}$ is $V^{\star}_{\gamma}(x):=x^{\top}P^{\star}_{\gamma}x$, where $P^{\star}_{\gamma}$ is a symmetric matrix. Hence SA5 holds with $\overline{\alpha}_{V^{\star}}(s)=\lambda_{\max}(P^{\star}_{1}) s$ and $\gamma^{\star}=1-\tfrac{\lambda_{\min}(S_1)}{\lambda_{\max}(P^{\star}_{1})}$. Provided we take $K_0$ such that $\gamma_0>\gamma^{\star}$, all the conditions of Corollary 1 are verified and we conclude that the system \eqref{Set function i} verifies a global exponential stability property for any $\gamma \in (\gamma^{\star},\gamma_0)$ after a sufficient number of iterations whose estimate is given by \eqref{estimate i}.

\subsection{Nonholonomic Integrator}
    
Consider the nonholonomic integrator as in \cite[Example 2]{Grimm-et-al-tac2005}, that is, 
\begin{equation}
\begin{split}
    x_{1}^{+} & = x_1+u_1\\
    x_{2}^{+} & = x_2 + u_2\\
    x_{3}^{+} & = x_3 +x_1u_2 - x_2u_1,
\end{split}
\end{equation}
where $x=(x_1,x_2,x_3)\in \R^3$, $u=(u_1,u_2) \in \mathcal{U}(x)=\R^2$. Let $\sigma(x)= x_{1}^2+x_{2}^2+10\euclidnorm{x_3}$ and $\ell(x,u)= \sigma(x) + \euclidnorm{u}^2$ for any $x\in \R^3$ and $u\in \R^2$. 

Thanks to \cite{Postoyan-et-al-tac(optimal)}, SA1, SA3 and item (i) of SA5 are verified with $W=0$, $\underline{\alpha}_W=0$, $\alpha_W=\mathbb{I}$ and $\alpha_{V^{\star}}=\tfrac{22}{5}\mathbb{I}$. We assume that SA2 holds. As the conditions of Lemma 1 are satisfied by taking $h(x)=(\tfrac{1}{15}x_1,-x_2)$ for any $x\in\R^{n_x}$ with $\chi=\mathbb{I}$, $a=\tfrac{256}{225}$, $M=\tfrac{22}{3}$. As a consequence, SA4 is satisfied with $\gamma_0=\tfrac{225}{256}$ and $\alpha_V(\cdot,\gamma) = \tfrac{M}{1-a\gamma}\mathbb{I}$ for any $\gamma\in(0,\gamma_0)$. Moreover, using the condition given in Corollary 1, item (ii) of SA5 is ensured with $\gamma^{\star}=\tfrac{17}{22}$. As $\tfrac{17}{22} \approx 0.77 < 0.88 \approx \tfrac{225}{256}$, we have $\gamma^{\star}<\gamma_0$ such that all the conditions of Corollary 1 are verified. We conclude that the system \eqref{Set function i} verifies a global exponential stability property for any $\gamma \in (\gamma^{\star},\gamma_0)$ after a sufficient number of iterations whose estimate is given by $i^{\star}_{\gamma} = \left \lceil{\tfrac{\ln{\big(\tfrac{(330\gamma-255)(1-\tfrac{256}{225}\gamma)}{21296\gamma(1-\gamma)^2}\big)}}{\ln{\big(\gamma-\tfrac{110\gamma-85}{484(1-\gamma)}\big)}}}\right \rceil$. Thus, when $\gamma=0.86$, $i^{\star}_{\gamma}=20$, which means that the systems controlled by the policies generated by PI exhibit an exponential stability property after 20 iterations.

\section{Conclusion}
We have analyzed the stability of general nonlinear discrete-time systems controlled by sequences of inputs generated by PI for an infinite-horizon discounted cost. Inspired by \cite{Granzotto-arxiv22}, novel near-optimal bounds have been established for the discounted problem. These new bounds do not blow up when the discount factor tends to 1 contrary to, e.g., \cite{Bertsekas-book12(adp),Munos-03}. The novel near-optimality bounds were then exploited to also provide general conditions under which PI generates stabilizing policies after sufficiently many iterations. 

In future work, we plan to relax some of the made assumptions in particular SA3 and item (ii) of SA5, and to investigate the case where $\gamma$ is increased with the number of iterations as done in \cite{Lamperski_2020} for the LQ problem.

\section*{Appendix}

\subsection{Proof of Proposition 1}
Let $\gamma \in \big(\gamma^{\star},\gamma_0\big)$, $x\in \R^{n_x}$ and $v=f\big(x,h^{\star}_{\gamma}(x)\big)$ with $h^{\star}_{\gamma}(x)\in H^{\star}_{\gamma}(x)$. Since $\ell$ is non-negative and by (i) of SA5,
    \begin{equation}
    \label{Sandwich bound for V_star}
        \ell\big(x,h^{\star}_{\gamma}(x)\big)\leq V^{\star}_{\gamma}(x) \leq \overline{\alpha}_{V^\star}\big(\sigma(x)\big).
    \end{equation}
These inequalities will be useful in the following. On the other hand, we have by definition of $V^{\star}_{\gamma}$
\begin{equation}
    \label{Bellman equation star}
    V^{\star}_{\gamma}(x)=\ell(x,h^{\star}_{\gamma}(x))+\gamma V^{\star}_{\gamma}(v),
\end{equation}
 therefore
\begin{equation}
\label{Lemma 2.1}
    V^{\star}_{\gamma}(v)-V^{\star}_{\gamma}(x)=-\tfrac{1}{\gamma}\ell\big(x,h^{\star}_{\gamma}(x)\big)+\tfrac{1-\gamma}{\gamma}V^{\star}_{\gamma}(x),\\
\end{equation}
which gives, using \eqref{Sandwich bound for V_star},
\begin{equation}
    \label{ineq difference V_star}
     V^{\star}_{\gamma}(v)-V^{\star}_{\gamma}(x)\leq-\tfrac{1}{\gamma}\ell\big(x,h^{\star}_{\gamma}(x)\big)+\tfrac{1-\gamma}{\gamma}\overline{\alpha}_{V^{\star}}\big(\sigma(x)\big).
\end{equation}
We define $Y^{\star}_{\gamma}:=V^{\star}_{\gamma}+\tfrac{1}{\gamma}W$. In view of SA3 and item (i) of SA5, $Y^{\star}_{\gamma}(x)\leq \overline{\alpha}_{Y^{\star}}\big(\sigma(x)\big)$ with $\overline{\alpha}_{Y^{\star}}:=\overline{\alpha}_{V^{\star}}+\tfrac{1}{\gamma^{\star}}\overline{\alpha}_W \in \mathcal{K}_{\infty}$. From the second inequality in \eqref{Detectability assumption} and the facts that $W(x)\ge0$ and $\gamma<1$, we have $\tfrac{1}{\gamma}W(x) \ge W(x) \ge \alpha_W\big(\sigma(x)\big)-\ell\big(x,h^{\star}_{\gamma}(x)\big)+W(v)\ge \alpha_W\big(\sigma(x)\big)-\ell\big(x,h^{\star}_{\gamma}(x)\big)$. Consequently, using \eqref{Sandwich bound for V_star}, $Y^{\star}_{\gamma}(x)\ge \alpha_W\big(\sigma(x)\big)-\ell\big(x,h^{\star}_{\gamma}(x)\big) + \ell\big(x,h^{\star}_{\gamma}(x)\big)=\alpha_W\big(\sigma(x)\big)=:\underline{\alpha}_{Y^{\star}}(\sigma(x))$ with $\underline{\alpha}_{Y^{\star}} \in \mathcal{K}_{\infty}$. On the other hand, in view of \eqref{Detectability assumption}, 
\begin{equation}
    \label{diff W1}
    \tfrac{1}{\gamma}W(v)-\tfrac{1}{\gamma}W(x)\leq -\tfrac{1}{\gamma}\alpha_W\big(\sigma(x)\big)+\tfrac{1}{\gamma} \ell\big(x,h^{\star}_{\gamma}(x)\big).
\end{equation}
Using \eqref{ineq difference V_star} and \eqref{diff W1}, we finally have 
\begin{equation} 
\label{Stability equation}
    Y^{\star}_{\gamma}(v)-Y^{\star}_{\gamma}(x)  \leq \tfrac{1}{\gamma}\Big(-\alpha_W\big(\sigma(x)\big) + (1-\gamma)\overline{\alpha}_{V^{\star}}\big(\sigma(x)\big)\Big).
\end{equation} 
In view of item (ii) of SA5, we have the next equivalences
\begin{align}
    (1-\gamma^{\star})\overline{\alpha}_{V^{\star}}\big(\sigma(x)\big)&\leq\alpha_W\big(\sigma(x)\big)\nonumber\\
    (1-\gamma)(1-\gamma^{\star})\overline{\alpha}_{V^{\star}}\big(\sigma(x)\big)&\leq(1-\gamma)\alpha_W\big(\sigma(x)\big)\nonumber\\
    (1-\gamma)(1-\gamma^{\star})\overline{\alpha}_{V^{\star}}\big(\sigma(x)\big) \nonumber \\
    - (1-\gamma^{\star})\alpha_W\big(\sigma(x)\big)&\leq-(\gamma-\gamma^{\star})\alpha_W\big(\sigma(x)\big)\nonumber\\
    \label{stability key inequality}
    (1-\gamma)\overline{\alpha}_{V^{\star}}\big(\sigma(x)\big) - \alpha_W\big(\sigma(x)\big) &\leq-\tfrac{\gamma-\gamma^{\star}}{1-\gamma^{\star}}\alpha_W\big(\sigma(x)\big).
\end{align}

As $\gamma \in (\gamma^{\star},\gamma_0)$, we have that $\alpha_{Y^{\star}}(\cdot,\gamma):=\tfrac{\gamma-\gamma^{\star}}{1-\gamma^{\star}}\alpha_W(\cdot) \in \mathcal{K_{\infty}}$. Moreover, as $Y^{\star}_{\gamma}(x)\leq \overline{\alpha}_{Y^{\star}}\big(\sigma(x)\big)$ with $\overline{\alpha}_{Y^{\star}} \in \mathcal{K}_{\infty}$, we derive from \eqref{Stability equation} and \eqref{stability key inequality},
\begin{equation}
    Y^{\star}_{\gamma}(v) \leq Y^{\star}_{\gamma}(x) - \tfrac{1}{\gamma} \alpha_{Y^{\star}}\Big(\overline{\alpha}^{-1}_{Y^{\star}}\big(Y^{\star}_{\gamma}(x)\big),\gamma\Big).
\end{equation}
Denoting $\widetilde{\alpha}_{Y^{\star}}(\cdot,\gamma):=\alpha_{Y^{\star}}(\overline{\alpha}_{Y^{\star}}^{-1}(\cdot),\gamma)$ and repeating the same reasoning, we have by induction that for any $k\in \Z_{\ge 0}$,
\begin{equation}
\label{Inductive composition}
    Y^{\star}_{\gamma}(\phi(k+1,x,h^{\star}_{\gamma}))\leq (\mathbb{I}-\tfrac{1}{\gamma}\widetilde{\alpha}_{Y^{\star}}(\cdot,\gamma))(Y^{\star}_{\gamma}(\phi(k,x,h^{\star}_{\gamma}))).
\end{equation}
 We can then apply \cite[Lemma 14]{Granzotto-arxiv22} with $\alpha(\cdot,\gamma)=\mathbb{I}-\tfrac{1}{\gamma}\widetilde{\alpha}_{Y^{\star}}(\cdot,\gamma)$ and $s_k(x)=Y^{\star}_{\gamma}(\phi(k,x,h^{\star}_{\gamma}))$ for any $k\in \Z_{\ge0}$ and $x\in \R^{n_x}$. We note using \eqref{Inductive composition} that $s_{k+1}(x)\leq\alpha(s_k(x),\gamma)$ as required. Moreover as $\mathbb{I}-\tfrac{1}{\gamma}\widetilde{\alpha}_{Y^{\star}}(\cdot,\gamma)$ is continuous, zero at zero, nonnegative and $\alpha(s,\gamma)<s$ for all $s >0$. There exists $\widetilde{\beta}^{\star}_{\gamma}\in \mathcal{KL}$ independant of $x$ and $v$ such that
 \begin{equation}
    \label{Lemma 14 MG}
     Y^{\star}_{\gamma}(\phi(k,x,h^{\star}_{\gamma}))\leq \widetilde{\beta}^{\star}_{\gamma}(Y^{\star}_{\gamma}(x),k),
 \end{equation}
 with $\widetilde{\beta}^{\star}_{\gamma} : (s,k)\mapsto \max_{\widehat{s}\in [0,s]}\alpha^{(k)}(\widehat{s},\gamma)$. As $\underline{\alpha}_{Y^{\star}}\big(\sigma(x)\big)\leq Y^{\star}_{\gamma}(x)\leq \overline{\alpha}_{Y^{\star}}\big(\sigma(x)\big)$, by \eqref{Lemma 14 MG} we obtain
 \begin{align}
     \underline{\alpha}_{Y^{\star}}\big(\sigma(\phi^{\star}_{\gamma}(k,x))\big)&\leq \widetilde{\beta}^{\star}_{\gamma}\big(\overline{\alpha}_{Y^{\star}}(\sigma(x)),k\big)\nonumber\\
     \sigma\big(\phi^{\star}_{\gamma}(k,x)\big)&\leq \beta^{\star}_{\gamma}\big(\sigma(x),k\big),
 \end{align}
with $\beta^{\star}_{\gamma}(s,k) \mapsto \underline{\alpha}_{Y^{\star}}^{-1}(\widetilde{\beta}^{\star}_{\gamma}(\overline{\alpha}_{Y^{\star}}(s),k)) \in \mathcal{KL}$. This concludes the proof. \hfill $\blacksquare$

\subsection{Proof of Theorem 1}
    Let $x\in \R^{n_x}$, $i \in \Z_{>0}$, $h^{i}_{\gamma}\in H^{i}_{\gamma}$, $h^{\star}_{\gamma}\in H^{\star}_{\gamma}$ and $\gamma \in (\gamma^{\star},\gamma_0\big)$. By Bellman equation and definition of $V^{i}_{\gamma}(x)$,
    \begin{align*}
        V^{i}_{\gamma}(x)-V^{\star}_{\gamma}(x)= & \, \ell\big(x,h^{i}_{\gamma}(x)\big)+\gamma V^{i}_{\gamma}\big(f(x,h^{i}_{\gamma}(x))\big)\\
       & - \ell\big(x,h^{\star}_{\gamma}(x)\big) -\gamma V^{\star}_{\gamma}\big(f(x,h^{\star}_{\gamma}(x))\big).
    \end{align*}
    As $h^{i}_{\gamma}(x)\in H^{i}_{\gamma}(x)=\underset{u\in \mathcal{U}(x)}{\text{argmin}}\{\ell(x,u)+\gamma V^{i-1}_{\gamma}\big(f(x,u)\big)\}$, we have for any $u\in \mathcal{U}(x)$, using Lemma 2,
    \begin{align*}
        V^{i}_{\gamma}(x)-V^{\star}_{\gamma}(x)\leq & \,  \ell(x,u)+\gamma V^{i-1}_{\gamma}\big(f(x,u)\big)- \ell\big(x,h^{\star}_{\gamma}(x)\big)\\
        &-\gamma V^{\star}_{\gamma}\big(f(x,h^{\star}_{\gamma}(x))\big).
    \end{align*}
     Thus by taking $u=h^{\star}_{\gamma}(x)\in H^{\star}_{\gamma}(x)$, we have
    \begin{equation}
        V^{i}_{\gamma}(x)-V^{\star}_{\gamma}(x)\leq \gamma \big(V^{i-1}_{\gamma}-V^{\star}_{\gamma}\big)(\phi(1,x,h^{\star}_{\gamma})).
    \end{equation}
    Repeating the above reasoning $i-1$ times, we obtain
    \begin{equation}
        \label{Prop 2-1}
        V^{i}_{\gamma}(x)-V^{\star}_{\gamma}(x)\leq  \gamma^i \big(V^{0}_{\gamma}-V^{\star}_{\gamma}\big)(\phi(i,x,h^{\star}_{\gamma})).
    \end{equation}
We have shown that the first inequality of \eqref{Convergence PI 1} holds. We now prove the second inequality of \eqref{Convergence PI 1}. Since $V^{\star}_{\gamma}\ge 0$, by SA4, we have 
\begin{equation}
\label{Ineq diff V1}
    V^{0}_{\gamma}(x)-V^{\star}_{\gamma}(x)\leq V^{0}_{\gamma}(x) \leq \overline{\alpha}_{V}\big(\sigma(x),\gamma\big).
\end{equation}
Thus, using \eqref{Ineq diff V1}, \eqref{Prop 2-1}, Proposition 1 and the fact that $\overline{\alpha}_{V}(\cdot,\gamma)$ is non-decreasing, we finally have 
\begin{align*}
    V^{i}_{\gamma}(x)-V^{\star}_{\gamma}(x)& \leq \gamma^i\big(V^{0}_{\gamma}-V^{\star}_{\gamma}\big)(\phi(i,x,h^{\star}_{\gamma}))\\
    &\leq \gamma^i \overline{\alpha}_{V}\big(\beta^{\star}_{\gamma}(\sigma(x),i), \gamma\big).
\end{align*}
This concludes the proof. \hfill $\blacksquare$

\subsection{Proof of Theorem 2}
   Let $i \in \Z_{\ge0}$, $\gamma \in \big(\gamma^{\star},\gamma_0\big)$, $x\in \R^{n_x}$, $v=f\big(x,h^{i}_{\gamma}(x)\big)$ with $h^{i}_{\gamma}(x)\in H^{i}_{\gamma}(x)$. By following similar arguments as in the proof of Proposition 1, item (i) of Theorem 2 holds with $Y^{i}_{\gamma}=V^{i}_{\gamma}+\tfrac{1}{\gamma}W$, $\underline{\alpha}_Y=\alpha_W$ and $\overline{\alpha}_{\gamma}(\cdot,\gamma)=\overline{\alpha}_V(\cdot,\gamma)+\tfrac{1}{\gamma}\overline{\alpha}_W(\cdot)$.

On the other hand, a similar inequality as \eqref{Lemma 2.1} is verified,
\begin{equation}
    \label{Th 1-0}
    V^{i}_{\gamma}(v)-V^{i}_{\gamma}(x)=-\tfrac{1}{\gamma}\ell\big(x,h^{i}_{\gamma}(x)\big)+\tfrac{1-\gamma}{\gamma}V^{i}_{\gamma}(x).
\end{equation}
In view of \eqref{Detectability assumption}, 
\begin{equation}
    \label{diff W}
    \tfrac{1}{\gamma}W(v)-\tfrac{1}{\gamma}W(x)\leq -\tfrac{1}{\gamma}\alpha_W\big(\sigma(x)\big)+\tfrac{1}{\gamma} \ell\big(x,h^{i}_{\gamma}(x)\big).
\end{equation}
Using \eqref{Th  1-0} and \eqref{diff W}, we have 
\begin{align}
    Y^{i}_{\gamma}(v)-Y^{i}_{\gamma}(x)&= V^{i}_{\gamma}(v)-V^{i}_{\gamma}(x)+\tfrac{1}{\gamma}W(v) -\tfrac{1}{\gamma} W(x) \nonumber \\
     & \leq \tfrac{1-\gamma}{\gamma}\big(V^{i}_{\gamma}(x)-V^{\star}_{\gamma}(x)\big) -\tfrac{1}{\gamma} \alpha_W\big(\sigma(x)\big) \nonumber\\
     & \quad + \tfrac{1-\gamma}{\gamma}V^{\star}_{\gamma}(x).
\end{align}
Using Theorem 1 and item (i) of SA5,
\begin{align}
\label{Th 2 - 1}
    Y^{i}_{\gamma}(v)-Y^{i}_{\gamma}(x)& \leq \tfrac{1}{\gamma}\big((1-\gamma)\gamma^{i}\overline{\alpha}_{V}\big(\beta^{\star}_{\gamma}(\sigma(x),i),\gamma\big) \nonumber \\
     & \quad - \alpha_W(\sigma(x)) + (1-\gamma)\overline{\alpha}_{V^{\star}}(\sigma(x))\big).
\end{align}
In view of  \eqref{Assumption 3}, \eqref{stability key inequality} and recalling that $\alpha_{Y^{\star}}(\cdot,\gamma)=\tfrac{\gamma-\gamma^{\star}}{1-\gamma^{\star}}\alpha_W(\cdot)$,
\begin{align}
\label{Th 2-12}
    Y^{i}_{\gamma}(v)-Y^{i}_{\gamma}(x) & \leq \tfrac{1}{\gamma}\big((1-\gamma)\gamma^{i}\overline{\alpha}_{V}\big(\beta^{\star}_{\gamma}(\sigma(x),i)\big) \nonumber\\
    & \quad -\alpha_{Y^{\star}}(\sigma(x),\gamma)\big) \nonumber\\
    & \leq \tfrac{1}{\gamma} \big( -\alpha_Y(\sigma(x),\gamma)+\Upsilon^{i}(\sigma(x),\gamma)\big),
\end{align}
with $\alpha_Y(\cdot,\gamma):=\alpha_{Y^{\star}}(\cdot,\gamma)\in \mathcal{K}_{\infty}$ and $\Upsilon^i(s_1,s_2):=(1-s_2)s_2^i\overline{\alpha}_{V}\big(\beta^{\star}_{\gamma}(s_1,i)\big)$ for any $s_1 \in \R_{\ge0}$ and $s_2\in (\gamma^{\star},\gamma_0)$. We have obtained item (ii) of Theorem 2. \hfill $\blacksquare$

\subsection{Proof of Theorem 3}
    Let $\delta, \Delta>0$, $\gamma\in (\gamma^{\star},\gamma_0)$,  $x\in \R^{n_x}$ such that $\sigma(x)\leq \Delta$, $i\ge i^{\star}_{\gamma}$, $v=f\big(x,h^{i}_{\gamma}(x)\big)$ where $h^{i}_{\gamma}\in H^{i}_{\gamma}$. We define $\widetilde{\Delta}_{\gamma}:=\overline{\alpha}_Y(\Delta,\gamma)>0$ and $\widetilde{\delta}:=\underline{\alpha}_Y(\delta)>0$ where $\underline{\alpha}_Y$ and $\overline{\alpha}_Y$ come from Theorem 2. Using item (ii) of Theorem 2, we have
\begin{equation}
\label{Th 2-13}
    Y^{i}_{\gamma}(v)-Y^{i}_{\gamma}(x)\leq \tfrac{1}{\gamma} \big( -\alpha_Y(\sigma(x),\gamma)+\Upsilon^{i}(\sigma(x),\gamma)\Big).
\end{equation}
As $\Upsilon^{i}(\cdot,\gamma)$ is non-decreasing and $\alpha_Y(\cdot,\gamma) \in \mathcal{K}_{\infty}$, using item (i) of theorem 2 and the fact that $Y^{i}_{\gamma}(x)\leq \widetilde{\Delta}_{\gamma}$, \eqref{Th 2-13} leads to 
\begin{align}
    Y^{i}_{\gamma}(v)-Y^{i}_{\gamma}(x) &\leq \tfrac{1}{\gamma}\big(\Upsilon^i(\underline{\alpha}^{-1}_{Y}(Y^i_{\gamma}(x)),\gamma)\nonumber \\
    & \qquad- \alpha_Y(\overline{\alpha}_{Y}^{-1}(Y^{i}_{\gamma}(x),\gamma),\gamma)\big)\nonumber\\
    \label{Th 2-2}
    Y^{i}_{\gamma}(v)-Y^{i}_{\gamma}(x) &\leq \tfrac{1}{\gamma}\big(\Upsilon^i(\underline{\alpha}^{-1}_{Y}(\widetilde{\Delta}_{\gamma}),\gamma) \nonumber\\
    & \qquad- \alpha_Y(\overline{\alpha}_{Y}^{-1}(Y^{i}_{\gamma}(x),\gamma),\gamma)\big).
\end{align}
On the other hand, as $\beta^{\star}_{\gamma}\in \mathcal{KL}$, for any $s\in \R_{\ge0}$ and $i \in \Z_{\ge0}$ we have
\begin{equation}
    \beta^{\star}_{\gamma}(s,i)\leq \beta^{\star}_{\gamma}(s,0).
\end{equation}
Hence, as $\overline{\alpha}_{V}$, $\alpha_Y \in \mathcal{K}_{\infty}$ and as we chose $i^{\star}_{\gamma}\ge\tfrac{\ln{\big(\tfrac{\alpha_{Y}(\overline{\alpha}_{Y}^{-1}(\underline{\alpha}_{Y}(\delta),\gamma),\gamma)}{2(1-\gamma)\overline{\alpha}_{V}(\beta^{\star}_{\gamma}(\underline{\alpha}_{Y}^{-1}(\overline{\alpha}_{Y}(\Delta,\gamma)),0),\gamma)}}\big)}{\ln{(\gamma)}}$ , we have $\Upsilon^{i^{\star}_{\gamma}}\big(\underline{\alpha}^{-1}_{Y}(\widetilde{\Delta}_{\gamma}),\gamma\big) \leq (1-\gamma)\gamma^{i^{\star}_{\gamma}}\overline{\alpha}_V(\beta^{\star}_{\gamma}(\underline{\alpha}_{Y}^{-1}(\widetilde{\Delta}_{\gamma}),0),\gamma) \leq \tfrac{1}{2} \alpha_Y\big(\overline{\alpha}_{Y}^{-1}(\widetilde{\delta},\gamma),\gamma\big)$. Consequently, for any $i\ge i^{\star}_{\gamma}$, as $ \Upsilon^{i}\big(\underline{\alpha}^{-1}_{Y}(\widetilde{\Delta}_{\gamma}),\gamma\big) \leq \Upsilon^{i^{\star}_{\gamma}}\big(\underline{\alpha}^{-1}_{Y}(\widetilde{\Delta}_{\gamma}),\gamma\big)$ we have when $Y^{i}_{\gamma}(x)\ge\widetilde{\delta}$,
\begin{align}
    \Upsilon^{i}\big(\underline{\alpha}^{-1}_{Y}(\widetilde{\Delta}_{\gamma}),\gamma\big)& \leq \tfrac{1}{2} \alpha_Y\big(\overline{\alpha}_{Y}^{-1}(\widetilde{\delta},\gamma),\gamma\big)\nonumber\\
    & \leq \tfrac{1}{2}\alpha_Y\big(\overline{\alpha}_{Y}^{-1}(Y^{i}_{\gamma}(x),\gamma),\gamma\big).
\end{align}
We then derive from  \eqref{Th 2-2} that when $Y^{i}_{\gamma}(x)\ge\widetilde{\delta}$,
\begin{equation}
    Y^{i}_{\gamma}(v)-Y^{i}_{\gamma}(x) \leq -\tfrac{1}{2\gamma}\alpha_Y\big(\overline{\alpha}_{Y}^{-1}(Y^{i}_{\gamma}(x),\gamma),\gamma\big).
\end{equation}
As $\tfrac{1}{2\gamma}\alpha_Y\big(\overline{\alpha}_{Y}^{-1}(\cdot,\gamma),\gamma\big)\in \mathcal{K}_{\infty}$, by proceeding by iteration and using \cite[Theorem 8]{Nesic_Teel_Sontag_99} we deduce that there exists $\widetilde{\beta}_{\gamma}\in \mathcal{KL}$ such that for any solution $\phi$ to \eqref{Set function i} initialized at $x$ and any $k\in \Z_{\ge 0}$,
\begin{equation}
    Y^{i}_{\gamma}\big(\phi(k,x,h^{i}_{\gamma})\big)\leq \max\{\widetilde{\beta}_{\gamma}\big(Y^{i}_{\gamma}(x),k\big),\widetilde{\delta}\}.
\end{equation}
Using item $(i)$ of Theorem 2 and the definition of $\widetilde{\delta}$, we obtain 
\begin{equation}
    \sigma\big(\phi(k,x,h^{i}_{\gamma})\big)\leq  \max\{\beta_{\gamma}\big(\sigma(x),k\big),\delta\},
\end{equation}
with $\beta_{\gamma}(s,k):=\underline{\alpha}_Y^{-1}\big({\widetilde{\beta}_{\gamma}}(\overline{\alpha}_Y(s,\gamma),k)\big)$. We have obtained the desired result as $h^{i}_{\gamma}$ has been arbitrarily selected in $H^{i}_{\gamma}$.

$ $\hfill $\blacksquare$

\subsection{Proof of Corollary 1}
      Let $\gamma \in (\gamma^{\star},\gamma_0)$, $x\in \R^{n_x}$, $i\ge i^{\star}_{\gamma}$, $v\in F^{i}_{\gamma}(x)$. Using Theorem 3, in particular in view of \eqref{Th 2 - 1},
    \begin{align}
        Y^{i}_{\gamma}(v)-Y^{i}_{\gamma}(x) \leq & \tfrac{1}{\gamma}\Big((1-\gamma)\gamma^{i}\overline{\alpha}_{V}\Big(\beta^{\star}_{\gamma}\big(\sigma(x),i\big),\gamma\Big)\nonumber\\ 
        & - \alpha_W\big(\sigma(x)\big) + (1-\gamma)\overline{\alpha}_{V^{\star}}\big(\sigma(x)\big)\Big).
    \end{align}
   Using similar lines as in \cite[Corollary 2]{Postoyan-et-al-tac(optimal)} and inequality \eqref{stability key inequality}, we prove that Proposition 1 is verified in this case with $\beta^{\star}_{\gamma}(s,k)=Ke^{-\lambda(\gamma) k} s$ where  $K=\tfrac{\gamma^{\star}\overline{a}_{V^{\star}}+\overline{a}_{W}}{\gamma^{\star} a_W}>0$ and $\lambda(\gamma)=-\ln{\Big(1-\tfrac{a_W(\gamma-\gamma^{\star})}{\gamma(1-\gamma)(\overline{a}_{V^{\star}}+\tfrac{1}{\gamma^{\star}}\overline{a}_W)}\Big)}>0$. By the conditions of Corollary 1 and the fact that $\overline{\alpha}_{V}(s,\gamma)=\overline{a}_V(\gamma)s$ for any $s>0$, this implies
   \begin{align}
       Y^{i}_{\gamma}(v)-Y^{i}_{\gamma}(x) &\leq \tfrac{1}{\gamma}\Big(K \overline{a}_V(\gamma) \big(1-\gamma)(\gamma e^{-\lambda(\gamma)})^i - a_W \nonumber\\
       & \qquad +(1-\gamma)\overline{a}_{V^{\star}}\Big)\sigma(x).
   \end{align}
As we have that $\gamma^{\star}=\tfrac{\overline{a}_{V^{\star}}-a_W}{\overline{a}_{V^{\star}}}$, item (ii) of SA5 is verified and as $\gamma\in(\gamma^{\star},\gamma_0)$ using \eqref{stability key inequality} we can write 
\begin{align}
       Y^{i}_{\gamma}(v)-Y^{i}_{\gamma}(x) &\leq \tfrac{1}{\gamma}\big(K \overline{a}_V(\gamma)(1-\gamma)(\gamma e^{-\lambda(\gamma)})^i \nonumber \\
       \label{equation 49}
       & \quad - \tfrac{a_W(\gamma-\gamma^{\star})}{\gamma(1-\gamma)}\big)\sigma(x).
   \end{align}
As we chose $i^{\star}_{\gamma}\ge\tfrac{\ln{\big(\tfrac{\gamma^{\star}(\gamma-\gamma^{\star})a_W^2}{2\gamma(1-\gamma)^2\overline{a}_V(\gamma)(\gamma^{\star}\overline{a}_{V^{\star}}+\overline{a}_W)}\big)}}{\ln{\big(\gamma - \tfrac{\gamma^{\star}(\gamma-\gamma^{\star})a_W}{(1-\gamma)(\gamma^{\star}\overline{a}_{V^{\star}}+\overline{a}_W)
    }\big)}}$, we have $K \overline{a}_V(\gamma) \big(1-\gamma)(\gamma e^{-\lambda(\gamma)})^i\leq \tfrac{a_W(\gamma-\gamma^{\star})}{2\gamma(1-\gamma)}$ for any $i\ge i^{\star}_{\gamma}$. Hence, \eqref{equation 49} becomes,
\begin{equation}
    Y^{i}_{\gamma}(v)-Y^{i}_{\gamma}(x) \leq -\tfrac{a_W(\gamma-\gamma^{\star})}{2\gamma(1-\gamma)}\sigma(x).
\end{equation}

Moreover, by Theorem 2 and the conditions of Corollary 1, we have $ Y^{i}_{\gamma}(x) \leq \big(\overline{a}_V(\gamma)+\tfrac{1}{\gamma}\overline{a}_W\big)\sigma(x)=: \overline{a}_Y(\gamma)\sigma(x)$ with $\overline{a}_Y(\gamma):=\overline{a}_V(\gamma)+\tfrac{1}{\gamma}\overline{a}_W>0$. We then write 
\begin{equation}
    Y^{i}_{\gamma}(v) \leq \big(1-\tfrac{a_W(\gamma-\gamma^{\star})}{2\gamma(1-\gamma)\overline{a}_Y(\gamma)}\big)Y^{i}_{\gamma}(x),
\end{equation}
with $1-\tfrac{a_W(\gamma-\gamma^{\star})}{2\gamma(1-\gamma)\overline{a}_Y(\gamma)}\in (0,1)$. Let $x\in \R^{n_x}$ and denote $\phi^{i}_{\gamma}(k,x)$ a corresponding solution to \eqref{Set function i} at time $k\in \Z_{\ge 0}$, it holds that $Y^{i}_{\gamma}\big(\phi^{i}_{\gamma}(k,x)\big)\leq \big(1-\tfrac{a_W(\gamma-\gamma^{\star})}{2\gamma(1-\gamma)\overline{a}_Y(\gamma)}\big)^k Y^{i}_{\gamma}(x)$. Using the fact that $Y^{i}_{\gamma}(x)\ge \alpha_W\big(\sigma(x)\big)=a_W\sigma(x)$ (in view of Theorem 2 and conditions of Corollary 1), we derive that Corollary 1 holds with $c_1(\gamma)=\tfrac{\overline{a}_Y(\gamma)}{a_W}=\tfrac{\gamma \overline{a}_V(\gamma)+\overline{a}_W}{\gamma a_W}$ and $c_2(\gamma)=-\ln{\big(1-\tfrac{a_W(\gamma-\gamma^{\star})}{2\gamma(1-\gamma)\overline{a}_Y(\gamma)}}\big)$. \hfill $\blacksquare$

\bibliography{bib_global} 
\end{document}